\documentclass{tran-l}

\newtheorem{theorem}{Theorem}[section]
\newtheorem{lemma}[theorem]{Lemma}
\newtheorem{prop}[theorem]{Proposition}
\newtheorem{cor}[theorem]{Corollary}

\theoremstyle{definition}
\newtheorem{definition}[theorem]{Definition}
\newtheorem{example}[theorem]{Example}

\theoremstyle{remark}
\newtheorem{remark}[theorem]{Remark}

\numberwithin{equation}{section}

\begin{document}
\title[Uniqueness of blow-up solutions]
{Boundary blow-up in nonlinear elliptic equations of
Bieberbach--Rademacher type}
\author[F.-C. C\^{\i}rstea]{Florica-Corina C\^{\i}rstea}
\address{Department of Mathematics, The Australian National University,
Canberra, ACT 0200, Australia}
\email{Florica.Cirstea@maths.anu.edu.au}
\thanks{The research of the first author was carried out at
Victoria University (Melbourne) with the support of the Australian Government through DETYA}
\author[V. R\u adulescu]{Vicen\c tiu R\u adulescu}
\address{Department of Mathematics, University of Craiova,
200585 Craiova, Romania} \email{radulescu@inf.ucv.ro}
\urladdr{http://inf.ucv.ro/\textasciitilde radulescu}
\begin{abstract}
We establish the uniqueness of the positive solution for equations
of the form $-\Delta u=au-b(x)f(u)$ in $\Omega$,
$u|_{\partial\Omega}=\infty$. The special feature is to consider
nonlinearities $f$ whose variation at infinity is \emph{not
regular} (e.g., $\exp(u)-1$, $\sinh(u)$, $\cosh(u)-1$,
$\exp(u)\log(u+1)$, $u^\beta \exp(u^\gamma)$, $\beta\in {\mathbb
R}$, $\gamma>0$ or $\exp(\exp(u))-e$) and functions $b\geq 0$ in
$\Omega$ vanishing on $\partial\Omega$. The main innovation
consists of using Karamata's theory not only in the
statement/proof of the main result but also to link the
non-regular variation of $f$ at
infinity with the blow-up rate of the solution near $\partial\Omega$. 
\end{abstract}
\maketitle

\section{Introduction}

Let $\Omega\subset {\mathbb R}^N$ $(N\geq 3)$ be a smooth bounded domain.
We consider semilinear elliptic problems under the following form
\begin{equation} \label{no1}
\Delta u=g(x,u) \quad \mbox{in}\ \Omega, \end{equation}
subject to the singular boundary condition
\begin{equation} \label{no2}
u(x)\to \infty\ \ \mbox{as}\ d(x):={\rm dist}\,(x,\partial\Omega)\to 0 \ \
(\mbox{in short, }u=\infty\ \mbox{on}\ \partial\Omega).
\end{equation}

The nonnegative solutions of (\ref{no1})+(\ref{no2}) are called {\em large}
(or {\em blow-up}) {\em solutions}.

The study of large solutions has been initiated in 1916 by Bieberbach \cite{bi}
for the particular case $g(x,u)=\exp(u)$ and $N=2$.
He showed that there exists a unique solution of (\ref{no1}) such that
$u(x)-\log (d(x)^{-2})$ is bounded as $x\to\partial\Omega$.
Problems of this type arise in Riemannian geometry; if a Riemannian metric of the form
$|ds|^2=\exp(2u(x))|dx|^2$ has constant Gaussian curvature $-c^2$ then
$\Delta u=c^2 \exp(2u)$. Motivated by a problem in mathematical physics,
Rademacher \cite{rad} continued the study of Bieberbach on smooth
bounded domains in $ {\mathbb R}^3 $.
Lazer--McKenna \cite{lm}
extended the results of Bieberbach and Rademacher for bounded domains
in ${\mathbb R}^N$ satisfying a uniformal external sphere condition and
for nonlinearities $g(x,u)=b(x)\exp(u)$, where $b$ is continuous and strictly
positive on $\overline{\Omega}$.

The interest in large solutions extended to $N$-dimensional
domains and for other classes of nonlinearities (see e.g.,
\cite{be}, \cite{bm}, \cite{cr}, \cite{ccm}, \cite{cras2004},
\cite{dh}, \cite{gls}, \cite{ke}, \cite{lm2}--\cite{os}).

Let $g(x,u)=f(u)$ where $f$ satisfies
\begin{equation} \tag{$A$}
f\in C^1[0,\infty),\ f'(s)\geq 0\ \mbox{for } s\geq 0,\ f(0)=0
\ \mbox{and}\  f(s)>0\ \mbox{for}\ s>0.
\end{equation}
In this case, Keller \cite{ke} and Osserman
\cite{os} proved that large solutions of (\ref{no1}) exist if and only if
\begin{equation} \tag{$A_0$}
\int_1^\infty  \frac{dt}{\sqrt{F(t)}} <\infty,\quad \mbox{where}\
F(t)=\int_0^t f(s)\,ds.
\end{equation}
In a celebrated paper, Loewner and Nirenberg \cite{ln}
linked the uniqueness of the blow-up solution to
the growth rate at the boundary.
Motivated by certain geometric problems, they
established the uniqueness for the case $f(u)=u^{\frac{N+2}{N-2}}$ $(N>2)$.
Bandle and Marcus \cite{bm} give results on asymptotic behaviour and uniqueness
of the large solution for more general nonlinearities including
$f(u)=u^p$ for any $p>1$. Theorem 2.3 in \cite{bm} proves that when $(A)$ holds
and
\begin{equation} \tag{$B$}
\exists \mu>0\ \mbox{and}\ s_0\geq 1\ \mbox{such that}\
f(\tau s)\leq \tau^{\mu+1} f(s)\ \forall \tau\in (0,1)\ \forall
s\geq s_0/\tau
\end{equation}
then for any large solution of $\Delta u=f(u)$ we have
\begin{equation} \label{mar}
\lim_{d(x)\to 0}\frac{u(x)}{Z(d(x))}=1 \end{equation}
where $Z$ is a chosen solution of
\begin{equation} \label{zed}
\left\{\begin{aligned}
& Z''(r)=f(Z(r)),\quad r\in (0,\delta)\ \mbox{for some}\ \delta>0\\
& Z(r)\to \infty\ \mbox{as}\ r\to 0^+.
\end{aligned}
\right.
\end{equation}
If, in addition,
$ f(\tau s)\leq \tau f(s)$, for all $\tau\in (0,1)$ and $s>0$,
then the uniqueness of large solutions takes place.
Lazer and McKenna \cite{lm2} consider the case when the $C^1$-function $f$
is either defined and positive on ${\mathbb R}$ or is defined
on $[a_0,\infty)$ with $f(a_0)=0$ and $f(s)>0$ for $s>a_0$.
They prove the uniqueness of large solutions to $\Delta u=f(u)$ in
$\Omega\subset {\mathbb R}^N$, $N>1$,
under the assumptions (see \cite[Theorem 3.1]{lm2}):
\begin{eqnarray}
&& \begin{aligned}
& \Omega\ \mbox{satisfies both
a uniform internal sphere condition
and a uniform} \\
& \mbox{external sphere condition with the same constant}\ R_1>0
\end{aligned}\\
&& f'(s)\geq 0\ \mbox{for } s \mbox{ in the domain of } f;\\
&& \mbox{there exists }
a_1\ \mbox{such that } f'(s)\ \mbox{is nondecreasing for } s\geq a_1;\label{ult}\\
&& \lim_{s\to \infty} f'(s)/\sqrt{F(s)}=\infty.
\end{eqnarray}
Moreover, the asymptotics
of the large solution is found in terms of a difference
\[ \lim_{d(x)\to 0} [u(x)-Z(d(x))]=0,\quad \mbox{for any}\ Z \
\mbox{satisfying}\ (\ref{zed}). \]

We are interested in large solutions
of (\ref{no1}) when $g(x,u)=b(x)f(u)-au $, i.e.,
\begin{equation}
\tag{$P$}
-\Delta u=au-b(x)f(u)\quad \mbox{in}\ \Omega,
\end{equation}
where $f\in C^1[0,\infty)$, $a\in {\mathbb R}$ and
$b\in C^{0,\mu}(\overline{\Omega})$ ($0<\mu<1$) satisfies $b\geq 0$,
$b\not\equiv 0$ in $\Omega$.

Many papers (see e.g., \cite{at}, \cite{da}--\cite{gls}) have been written about
Eq. ($P$), on a bounded domain or
${\mathbb R}^N$, when $f(u)=u^p$ ($p>1$). For this case of nonlinearity
and $b>0$ on $\overline{\Omega}$, Eq. ($P$) subject to
$u=0$ on $\partial\Omega$
is referred to as the logistic equation. It is known that it has a unique
positive solution if and only if $a>\lambda_1(\Omega)$, where
$\lambda_1(\Omega)$ is the first Dirichlet eigenvalue of $(-\Delta)$ in
$\Omega$. We mention that the logistic equation has been proposed as a model
for population density of a steady-state single species $u(x)$ when $\Omega$
is fully surrounded by inhospitable areas.
However, not until recently was the case of a
{\em degenerate} logistic type equation considered, which allows $b$ to vanish
on $\overline{\Omega}$ (see \cite{at}, \cite{ggll} and \cite{ddm}).
The understanding of the asymptotics
for positive solutions of the degenerate logistic equation
leads to the study of large solutions (we refer to \cite{ggll} and \cite{gls}).

Let $\Omega_0$
denote the interior of the zero set of $b$ in
$\Omega$, i.e.,
\[\Omega_0={\rm int}\,\{x\in \Omega:\ b(x)=0\}.\]

We assume throughout that $\Omega_0$ is connected,
$\partial\Omega_0$ satisfies the exterior cone condition
(possibly, $\Omega_0=\emptyset$), $\overline{\Omega}_0\subset
\Omega$ and $b>0$ on $\Omega\setminus \overline{\Omega}_0$. Note
that $b\geq 0$ on $\partial\Omega$.

Let $\lambda_{\infty,1}$ be the first Dirichlet eigenvalue of
$(-\Delta)$ in $\Omega_0$. Set $\lambda_{\infty,1}=\infty$ if
$\Omega_0=\emptyset$.

Alama and Tarantello \cite{at} find the maximal interval $I$ for the parameter $a$
such that ($P$), subject to $u=0$ on $\partial\Omega$, has a positive solution
$u_a$, provided that
\begin{equation} \tag{$A_1$}
f\geq 0\ \mbox{and}\ f(u)/u\ \mbox{is increasing on}\ (0,\infty).
\end{equation}
Moreover, for each $a\in I$, the solution $u_a$ is unique
(see \cite[Theorem A (bis)]{at}).

Theorem 1.1 in \cite{ccm} proves that if ($A_0$) and ($A_1$) are
fulfilled, then Eq. ($P$) has large solutions if and only if $a\in
(-\infty, \lambda_{\infty,1})$. The uniqueness and asymptotic
behaviour near $\partial\Omega$ prove to be very challenging in
the above generality.

In \cite{cr} we advance for the first time the idea of using the
regular variation theory arising in applied probability to study
the uniqueness of large solutions. There we consider the case when
$f'$ varies regularly at infinity (see Definition~\ref{d1}).

Note that there are many nonlinearities $f(u)$, such as
$\exp(u)-1$, $\sinh(u)$, $\exp(\exp(u))-e$, $\exp(u)\log(u+1)$,
which do not fall in the category treated by Theorem~1 in
\cite{cr}. Although some examples might fit into the framework of
\cite[Theorem~2.3]{bm} or \cite[Theorem~3.1]{lm2}, the uniqueness
and growth rate at the boundary for large solutions of ($P$) have
not yet been studied when $a\not=0$ and $b$ vanishes in $\Omega$
with $b\equiv 0$ on $\partial\Omega$.

Our purpose is to fill in this gap by analysing a wide range of
functions $f$ and $b$.
We develop the research line opened up in \cite{cr} to treat here
the case when $f$ does not vary regularly at infinity. Thus our
approach for the uniqueness is different from that of
Bandle--Marcus and Lazer--McKenna, being based on Karamata's
theory.

\section{Framework and main result}

We first recall some results from the Karamata regular variation
theory (see \cite{bgt}).

\begin{definition} \label{d1} A measurable function $R:[A,\infty)\to
(0,\infty)$, for some $A>0$, is called {\em regularly varying at
infinity of index} $\rho\in \mathbb{R}$, in short $R\in RV_\rho$,
provided that
\[\lim_{u\to \infty}\frac{R(\xi u)}{R(u)}=\xi^\rho,\qquad \forall \xi>0.\]
When the index $\rho$ is zero, we say that the function is {\em slowly varying}.
\end{definition}

From now on, we do not write at infinity when the regular variation occurs there.
Notice that the transformation $R(u)=u^\rho L(u)$ reduces regular variation
to slow variation.
Examples of slowly varying functions are given by:
\begin{itemize}
\item[(i)] Every measurable function on $[A,\infty)$ which has a positive
limit at $\infty$.
\item[(ii)] The logarithm $\log u$, its iterates $\log_m u$,
and powers of $\log_m u$.
\item[(iii)] ${\rm exp}\{(\log u)^{\alpha_1}(\log_2 u)^{\alpha_2}\ldots
(\log_m u)^{\alpha_m}\}$ where $\alpha_i\in (0,1)$
and ${\rm exp}\{\frac{\log u}{\log\log u}\}$.
\end{itemize}

\begin{prop}[Representation Theorem] \label{km}
The function $L(u)$ is slowly varying if and only if it can be written
in the form
\[L(u)=M(u){\rm exp}\,\left\{\int_B^u \frac{\phi(t)}{t}\,dt\right\},\quad
\forall u\geq B\]
for some $B>0$, where $\phi\in C[B,\infty)$ satisfies
$\lim_{u\to \infty}\phi(u)=0$
and $M(u) $ is measurable on $[B,\infty)$ such that
$\lim_{u\to \infty}M(u)=\widehat{M}\in (0,\infty)$.
\end{prop}
If $M(u)$ is replaced by $\widehat{M}$ then the new function, say
$\widehat{L}(u)$, is referred to as a {\em normalised}
slowly varying function.
We see that
$\phi(u)=\frac{u\widehat{L}'(u)}{\widehat{L}(u)}$, $\forall u\geq B$. Conversely,
any function $\widehat{L}\in C^1[B,\infty)$ which is positive and satisfies
$\lim_{u\to \infty}\frac{u\widehat{L}'(u)}{\widehat{L}(u)}=0$ is a normalised slowly
varying function.

Note that {\em any} slowly varying function $L(u)$ is asymptotic
equivalent to some normalised slowly varying function $\widehat{L}(u)$
(i.e., $\lim_{u\to \infty}L(u)/\widehat{L}(u)=1$).

The notion of regular variation can be extended to any real number.
For instance, we say that $R(u)$
is regularly varying (on the right)
at the origin with index $\rho\in {\mathbb R}$ (and write
$R\in RV_\rho(0+)$) if $R(1/u)\in RV_{-\rho}$.
Let $NRV_\rho(0+)$ (resp., $NRV_\rho$)
denote the set of all normalised regularly varying
functions at $0$ (resp., $\infty$) of index $\rho$.

By $f_1(x)\sim f_2(x)$ as $d(x)\to 0$ we mean that $\lim_{d(x)\to
0} \frac{f_1(x)}{f_2(x)}=1$.

Our main result is

\begin{theorem} \label{teo1}
Let {\rm ($A_1$)} hold and
$f\circ \mathfrak{L}\in RV_\rho$ $(\rho>0)$ for some
$\mathfrak{L}\in C^2[A,\infty)$ satisfying
$\lim_{u\to \infty}\mathfrak{L}(u)=\infty$ and
$\mathfrak{L}'\in NRV_{-1}$.
Suppose that
\begin{equation} \tag{$H$}
\begin{aligned}
& b(x)\sim \mathfrak{K}^2(d(x))\ \mbox{as}\
d(x)\to 0,\ \mbox{where}\
\mathfrak{K}\in NRV_\theta(0+),\ \mbox{for some } \theta\geq 0\ \mbox{and}\\
& \mathfrak{K}\ \mbox{is nondecreasing near the origin if }\theta=0.
\end{aligned}
\end{equation}

Then, for any $a<\lambda_{\infty,1}$, Eq. {\rm ($P$)} has a unique
large solution $u_a$. In addition,
the blow-up rate of $u_a$ at
$\partial\Omega$ can be expressed by
\begin{equation} \label{ii}
u_a(x) \sim (\mathfrak{L}\circ \Phi)(d(x))\ \ \mbox{as}\ d(x)\to
0,\quad \forall a<\lambda_{\infty,1}.
\end{equation}
The function $\Phi$ is defined as follows
\begin{equation}\label{rv2}
\int_{\Phi(t)}^\infty
\frac{ [\mathfrak{L}'(y)]^{ \frac{1}{2} }}
{ y^{\frac{\rho+1}{2}}[{\boldsymbol L}_f(y)]^{\frac{1}{2}}}\,dy=\int_0^t \mathfrak{K}(s)\,ds,\quad
\forall t\in (0,\beta)\ \mbox{with}\ \beta>0\ \mbox{small},
\end{equation}
where
${\boldsymbol L}_f$ is a
normalised
slowly varying function
such that
$\lim_{u\to \infty} \frac{f(\mathfrak{L}(u))}{u^\rho {\boldsymbol L}_f(u)}=1$.
\end{theorem}

Note that Theorem~\ref{teo1} brings a new insight into the
asymptotics of the large solution of ($P$) even in the case $a=0$
and $b=1$. For instance, the function which is used in (\ref{ii})
to estimate the blow-up rate of the solution near $\partial\Omega$
is {\em not chosen as a solution of} (\ref{zed}). This fact will
allow us, through Corollary~\ref{expo}, to illustrate the
explosion pattern followed by the large solution when the
nonlinearity $f$ is of the form (\ref{tyr}) at infinity and
satisfies $(A_1)$. In particular, if $\lim_{u\to
\infty}\frac{f(u)}{[\exp_m(u^{\frac{1}{\alpha}})]^\rho}=1$,
($\alpha,\rho>0$ and $m\geq 1$ an integer), then the unique large
solution of $\Delta u=f(u)$ satisfies $\frac{u(x)}{\Psi(d(x))}\to
1$ as $d(x)\to 0$, where
\[
\Psi(d(x))=\left\{ \begin{aligned} &
[\log(d(x)^{-\frac{2}{\rho}})]^\alpha, \ \ \mbox{if }
m=1,\\
& [\log_m (d(x)^{-1})]^\alpha, \ \ \mbox{if } m\geq 2.
\end{aligned} \right.\]
We set $\log_m (\cdot)=\underset{m\ {\rm
times}}{(\underbrace{\log\circ \ldots \circ \log})}(\cdot) $ and
$\exp_m (\cdot)=\underset{m\ {\rm times}}{(\underbrace{\exp\circ
\ldots \circ \exp})}(\cdot) $, ${\mathbb Z}\ni m\geq 1$. If
$f(u)=\exp_2(u)+\cos(\exp_2(u))$ for $u$ large and $(A_1)$ holds,
the uniqueness of large solutions for $\Delta u=f(u)$ cannot be
inferred from the Lazer--McKenna result, since condition
(\ref{ult}) fails. Nevertheless, the uniqueness is valid as we can
derive from either \cite[Theorem 2.3]{bm} or Theorem~\ref{teo1}.
But it is not transparent through (\ref{mar}) that the large
solution fulfills $\lim_{d(x)\to
0}\frac{u(x)}{\log_2\left(\frac{1}{d(x)}\right)}= 1$, as
Corollary~\ref{expo} proves.

\begin{remark} \label{vit}
We point out that
$\mathfrak{L}'\in NRV_{-1}$ with $\lim_{u\to \infty}
\mathfrak{L}(u)=\infty $ if and only if
\begin{equation} \label{how}
\mathfrak{L}(u)=C{\rm exp}\left\{\int_B^u
\frac{\ell(t)}{t}\,dt\right\},\quad
\forall u\geq B>0
\end{equation}
where $C>0$ is a constant and $\ell$ is a normalised slowly varying function
satisfying $\lim_{u\to \infty} \ell(u)=0$ and $\lim_{u\to \infty}\int_B^u
\frac{\ell(t)}{t}\,dt=\infty $.
Nontrivial examples of functions $\mathfrak{L}$
are: ${\rm exp}\{(\log u)^\gamma\}$, where $\gamma\in (0,1)$,
$ {\rm exp}\left\{\frac{\log u}{\log \log u}\right\}$,
and $(\log_m u)^\alpha$ with $\alpha>0$.
\end{remark}

The hypothesis $f\circ \mathfrak{L}\in RV_\rho$ ($\rho>0$) is equivalent to
the existence of $g\in RV_\rho$ so that $f(u)=g(\mathfrak{L}^\leftarrow(u))$,
for $u$ large (where $\mathfrak{L}^\leftarrow$
denotes the inverse of $\mathfrak{L}$).
By Proposition~0.8~(v),
$\mathfrak{L}^\leftarrow$ is {\em rapidly varying with index} $\infty$
($\mathfrak{L}^{\leftarrow}\in RV_\infty$), i.e.,
\[ \lim_{u\to \infty} \frac{\mathfrak{L}^{\leftarrow}(\lambda u)}
{\mathfrak{L}^{\leftarrow}(u)}=
\begin{cases}
0 & \text{if $\lambda\in (0,1)$},\\
1 & \text{if $\lambda=1$},\\
\infty & \text{if $\lambda>1$}.
\end{cases} \]
Therefore, for $g(u)=u^\rho$,
$f(u)=[ \mathfrak{L}^{\leftarrow}(u)]^\rho$
is rapidly varying with index $\infty$.

If $g\in NRV_\rho$, then
${\boldsymbol L}_f$ (which appears in (\ref{rv2})) can be taken as
$\frac{f(\mathfrak{L}(u))}{u^\rho}$. Moreover, $\frac{f(u)}{u}$
is increasing in a neighbourhood of infinity. For this, it is enough to see
that $\lim_{u\to \infty}\frac{uf'(u)}{f(u)}>1$. Indeed, using
(\ref{how}), we derive that
\[ \lim_{y\to \infty}
\frac{f'(\mathfrak{L}(y))\mathfrak{L}(y)}{f(\mathfrak{L}(y))}=
\lim_{y\to \infty} \frac{yg'(y)}{g(y)}\,\frac{\mathfrak{L}(y)}
{y\mathfrak{L}'(y)}=\rho\,\lim_{y\to \infty}\frac{\mathfrak{L}(y)}
{y\mathfrak{L}'(y)}=\infty.\]

Proposition~\ref{km} will provide
countless functions $g\in NRV_\rho$ and $\mathfrak{L}$ as in (\ref{how}).
Hence, by taking
$f(u)=g(\mathfrak{L}^\leftarrow(u))$ ($u\geq B>0 $), the assumptions
of Theorem~\ref{teo1} are fulfilled. It remains only to
extend the definition of $f$ to the remaining part of $(0,\infty)$
such that the smoothness of $f$ and $(A_1)$ hold.

Regarding the assumption $(H)$, $\mathfrak{K}\in NRV_\theta(0+)$
if and only if there exists a normalised
slowly varying function $L_{\mathfrak{K}}$
such that
\begin{equation} \label{key}
\mathfrak{K}(t)=t^\theta L_{\mathfrak{K}}(1/t),\quad t\in (0,\nu)\ \mbox{with}\
\nu>0.
\end{equation}
Therefore (\ref{key})
is equivalent to saying that for some constants
$c,d>0$ and $\varphi\in C(0,\nu)$
with $\lim_{t\to 0^+}\varphi(t)=0$ we have
\[ \mathfrak{K}(t)=ct^\theta \exp\left(\int_t^d \frac{\varphi(y)}{y}\,dy\right),\quad
\mbox{for } t\in (0,d).\]
Some examples of $\mathfrak{K}$ as in $(H)$ are: $t^\theta$,
$(\sin t)^\theta$, $t^\theta/\exp \left[\frac{\log(1/t)}{\log\log(1/t)}\right]$,
$t^\theta/\exp[(-\log t)^\gamma]$ with $\gamma\in (0,1)$,
$t^\theta [\log(t+1)]^\alpha$ or
$t^\theta[\log_m(1/t)]^{-\alpha}$ with $\alpha>0$ and $m\geq 1$
an integer.

\begin{remark} \label{mmm} If in Theorem~\ref{teo1} we replace
$f\circ \mathfrak{L}\in RV_\rho$ by the hypothesis
$f'\in RV_\rho$ ($\rho>0$), then ($P$) still has a unique large solution
$u_a$, $\forall a<\lambda_{\infty,1}$. However,
the blow-up
rate of $u_a$ near $\partial\Omega$ is as follows (see
\cite[Theorem~1]{cr})
\begin{equation} \label{cl}
\lim_{d(x)\to 0}\frac{u_a(x)}{\left(\frac{2\theta+\rho+2}{(2+\rho)(\theta+1)}
\right)^{1/\rho}h(d(x))}=1,\quad \forall a<\lambda_{\infty,1}
\end{equation}
where $h$ is defined by
\begin{equation} \label{er}
\int_{h(t)}^\infty \frac{ds}{\sqrt{2F(s)}}=\int_0^t
\mathfrak{K}(s)\,ds,\quad \forall t\in (0,\nu).
\end{equation}
\end{remark}

\begin{remark} \label{ve}
The variation of $f$ at $\infty$
is {\em not regular} in Theorem~\ref{teo1}
(i.e., $f\not\in RV_{\gamma}$, for any $\gamma\in {\mathbb R}$) in contrast
to Remark~\ref{mmm} where $f\in NRV_{\rho+1}$.
This fact will bring a significant change in the explosion speed of the
large solution of ($P$).

By Lemma~\ref{pro} we know that
$\Phi\in NRV_{\frac{-2(\theta+1)}{\rho}}(0+)$.
Since $\mathfrak{L}$ varies slowly at infinity, we
can invoke \cite[Proposition~0.8 (iv)]{res} to conclude that
$\mathfrak{L}\circ \Phi\in RV_0(0+)$.

We show that, in the setting of Remark~\ref{mmm},
$h\in RV_{\frac{-2(\theta+1)}{\rho}}(0+)$.
It is easy to check that
$T(u)=(\int_0^{1/u}\mathfrak{K}(s)\,ds)^{-1}\in RV_{\theta+1}$.
Set $Y(u)=(\int_u^\infty \frac{ds}{\sqrt{2F(s)}})^{-1}$, for $u>0$.
Clearly, $Y$ is increasing on $(0,\infty)$, $Y(\infty)=\infty$ and
$Y\in RV_{\rho/2}$. By (\ref{er}), we find
$h(1/u)=Y^{\leftarrow}(T(u))$ for $u$ sufficiently large, where
$Y^{\leftarrow}(u)$ is the inverse of $Y(u)$.
By Proposition~0.8(v) in \cite{res}, $Y^\leftarrow\in RV_{2/\rho}$ so that
$h(1/u)\in RV_{\frac{2(1+\theta)}{\rho}}$.
\end{remark}

As a consequence of Theorem~\ref{teo1} and (\ref{jar}), we obtain

\begin{cor} \label{expo}
Let {\rm ($A_1$)} and $(H)$ hold. Assume that there exists
$\alpha,\rho>0$ and an integer $m\geq 1$ such that $f((\log_m
u)^\alpha)\in RV_\rho$.

Then Eq. {\rm ($P$)} has a unique large solution $u_a$, for any
$a<\lambda_{\infty,1}$. Moreover,
\begin{equation} \label{son}
\lim_{d(x)\to 0}\frac{u_a(x)}{\left[\log_m\left(\frac{1}{d(x)}
\right)\right]^\alpha}=
\begin{cases}
\left(\frac{2(1+\theta)}{\rho}\right)^{\alpha} & \text{if $m=1$},\\
1 & \text{if $m\geq 2$}.
\end{cases}
\end{equation}
\end{cor}

\begin{remark}
For $m=1$ the influence of $f$ (resp., $\mathfrak{K}$) into the
blow-up rate (\ref{son}) of the large solution can be seen through
$\alpha$ and $\rho$ (resp., $\theta$). Nevertheless, if $m\geq 2$,
then the order of iteration for logarithm changes accordingly in
the asymptotic behaviour (\ref{son}) that proves to be {\em
independent of the index of regular variation} $\rho$ (for
$f((\log_m u)^\alpha)$) and $\theta$ (for $\mathfrak{K}$).
\end{remark}

The assumption $f((\log_m u)^\alpha)\in RV_\rho$ holds if and only
if there exists a slowly varying function $L$ such that
\begin{equation} \label{tyr}
f(u)=\left[\exp_m(u^{\frac{1}{\alpha}})\right]^{\rho}
L(\exp_m(u^{\frac{1}{\alpha}})),\quad u\geq B>0.
\end{equation}
Such examples are given below:
\begin{itemize}
\item[(i)] $f(u)=u^{\beta}{\rm exp}\{\rho u^{\frac{1}{\alpha}}\}$,
$f(u)=(\log u)^\beta {\rm exp}\{\rho u^{\frac{1}{\alpha}}\}$,
where $\beta\in {\mathbb R}$ is arbitrary; \item[(ii)] $f(u)={\rm
exp}\left\{ u^{\frac{1}{\alpha}}(\rho+\frac{\alpha}{\log u}
)\right\}$, $f(u)={\rm exp}\left\{ u^{\frac{1}{\alpha}}
[\rho+u^{\frac{-2}{3\alpha}}
\cos(u^{\frac{1}{3\alpha}})]\right\}$; \item[(iii)]
$f(u)=\exp\left\{u^{\frac{1}{\alpha}}(\rho+u^{\frac{\alpha_1-1}{\alpha}})
\right\}$ with $\alpha_1\in (0,1)$; \item[(iv)]
$f(u)=\exp\{u^{\frac{1}{\alpha}}+\rho
\exp\{u^{\frac{1}{\alpha}}\}\}$, $f(u)=\exp \{
(u^{-\frac{1}{\alpha}}+\rho) \exp\{u^{\frac{1}{\alpha}}\}\}$
\end{itemize}
($m=1$ in (i)--(iii) and $m=2$ in (iv)).

\begin{example} Among functions $f$ which fulfill the hypotheses of
Corollary~\ref{expo}, we illustrate: $f(u)=\exp\{u\}-1$,
$f(u)=\sinh(u)$, $f(u)=\cosh(u)-1$, $f(u)=\exp\{u\}\log(u+1)$,
$f(u)=u^\beta\exp\{\rho u^{\frac{1}{\alpha}}\}$ with $\beta\geq
1$, $\alpha,\rho>0$, $f(u)=u^\beta \exp(\exp\{u\})$ with
$\beta\geq 1$ and $f(u)=\exp(\exp\{u\})-e$.
\end{example}

Boundary blow-up phenomena for $(P)$ with $a=0$, $b=1$ and
$f(u)=u^p$, $1<p\leq 2$, appear in the analytical theory of a
Markov process called superdiffusion. In this case, the uniqueness
of the large solution was studied in Dynkin \cite{dyn1, dyn2} by
probabilistic techniques. It is remarkable that Dynkin's papers
realize, on one hand, a connection between superprocesses and
singularity phenomena and, on the other hand, they contain a
probabilistic representation of the minimal large solution. By
means of a probabilistic representation, a uniqueness result in
domains with non-smooth boundary was established by Le Gall
\cite{gall} in the case $p=2$. The existence of large solutions is
usually deduced by comparison methods combined with
Keller-Osserman {\it a~priori} bounds, Calderon-Zygmund estimates,
Agmon-Douglis-Nirenberg's theory, or Alexandrov and Krylov-Safonov
techniques.

Our interest falls here on the uniqueness of large solutions to
$(P)$ when $f$ does not vary regularly at infinity (thus excluding
the power case). Note that if $f(u)=\exp(u)-1$ or
$f(u)=\exp(u)-u-1$, then by Corollary~\ref{expo} the equation
$\Delta u=f(u)$ in $\Omega$ has a unique large solution which
satisfies $\lim_{d(x)\to 0}\frac{u(x)}{\log[d(x)]^{-2}}=1$. This
asymptotic behaviour is exactly the same as for the unique large
solution of $\Delta u=\exp(u)$ in $\Omega$, going back to the
pioneering works of Bieberbach \cite{bi} and Rademacher
\cite{rad}. For the two-term asymptotic expansion of the large
solution of $\Delta u=\exp(u)$ we refer to \cite{bandle}. We point
out that our approach is completely different from the above
papers for it relies exclusively on the regular variation theory
(see \cite{bgt} for details) not only in the statement, but also
in the proof of the main result.

\section{Auxiliary results}

For details about Propositions~\ref{reg} and \ref{ka} we refer the reader
to \cite{bgt} (\cite{res} or \cite{se}).

\begin{prop}[Elementary properties of slowly varying functions]
\label{reg} Assume that $L$ is a slowly varying function. Then the following
hold
\begin{itemize}
\item[(i)] $\log L(u)/\log u\to 0$ as $u\to \infty$.
\item[(ii)] For any $m>0$, $u^m L(u)\to \infty$, $u^{-m}L(u)\to 0$
as $u\to \infty$.
\item[(iii)] $(L(u))^m$ varies slowly for every $m\in \mathbb{R}$.
\item[(iv)] If $L_1$ varies slowly, so do $L(u)L_1(u)$ and $L(u)+L_1(u)$.
\end{itemize}
\end{prop}

\begin{remark} \label{ops}
If $g\in RV_\rho$ with $\rho>0$ ($\rho<0$), then
$\lim_{u\to \infty} g(u)=\infty$ ($0$).
However, the behaviour
at infinity for a slowly varying function cannot be predicted.
We see that $L(u)={\rm exp}\{(\log u)^{1/3}\cos((\log u)^{1/3})\}$
is a (normalised) slowly varying function (use
$\lim_{u\to \infty}\frac{uL'(u)}{L(u)}=0$) for which
$\liminf_{u\to \infty}L(u)=0$ and $\limsup_{u\to \infty} L(u)=\infty$.
\end{remark}

\begin{prop}[Karamata's Theorem] \label{ka}
Let $R\in RV_\rho$ be locally bounded in $[A,\infty)$.
Then, for any $j<-(\rho+1)$ (resp., $j=-(\rho+1)$ if
$\int^\infty x^{-(\rho+1)}R(x)\,dx<\infty $)
\[ \lim_{u\to \infty} \frac{u^{j+1}R(u)}{\int_u^\infty x^j R(x)\,dx}=
-(j+\rho+1). \]
\end{prop}

Under the assumptions of Theorem~\ref{teo1}, we prove

\begin{lemma} \label{pro}
The function $\Phi$ given by {\rm (\ref{rv2})} is well defined on some interval
$(0,\beta)$. Furthermore, $\Phi\in NRV_{\frac{-2(\theta+1)}{\rho}}(0+)$
satisfies
\begin{equation} \label{jar}
\lim_{t\to 0^+}\frac{\log_m \Phi(t)}{\log_m(\frac{1}{t})}=
\begin{cases}
\frac{2(1+\theta)}{\rho} & \text{if $m=1$},\\
1 & \text{if $m\geq 2$}.
\end{cases}
\end{equation}
\begin{equation}
\lim_{t\to 0^+}\frac{\Phi(t)\Phi''(t)}{[\Phi'(t)]^2}=1+
\frac{\rho}{2(\theta+1)} \quad \mbox{and} \quad
\lim_{t\to 0^+}
\frac{\mathfrak{L}(\Phi(t))}{\mathfrak{L}'(\Phi(t))}\,\frac{\Phi(t)}{[\Phi'(t)]^2}=0.\label{ch2}
\end{equation}
\end{lemma}

\begin{proof} Let $b>0$ be such that
${\boldsymbol L}_f$ resp., $\mathfrak{L}'$ is positive
on $[b,\infty)$. Since $\mathfrak{L}'\in RV_{-1}$ and
${\boldsymbol L}_f$ is slowly varying, Proposition~\ref{reg} yields
\[ \lim_{u\to \infty}
\frac{[\mathfrak{L}'(u)]^{\frac{1}{2}} }
{ u^{\frac{\rho+1}{2}} [{\boldsymbol L}_f(u)]^{\frac{1}{2}}}\, u^{1+\tau}=0,
\quad \mbox{for any}\ \tau\in (0,\rho/2).\]
Therefore, there exists $B>b$ large so that
\[ \zeta(x)=\int_x^\infty
\frac{[\mathfrak{L}'(y)]^{\frac{1}{2}} }
{ y^{\frac{\rho+1}{2}} [{\boldsymbol L}_f(y)]^{\frac{1}{2}}}
\,dy<\infty,\quad
\forall x>B.\]
It follows that
$\Phi$ is well defined on $(0,\beta)$, for some $\beta>0$.
Moreover, $\Phi\in C^2(0,\beta)$ and $\lim_{t\to 0^+}\Phi(t)=
\infty$.
Using (\ref{rv2}), we find
\begin{equation} \label{et1}
\frac{-\Phi'(t) [\mathfrak{L}'(\Phi(t))]^{\frac{1}{2}}}
{[\Phi(t)]^{\frac{\rho+1}{2}}[{\boldsymbol L}_f(\Phi(t))]^{\frac{1}{2}}}=\mathfrak{K}(t),\quad 
\forall
t\in (0,\beta).\end{equation}
In view of Proposition~\ref{ka}, we have
\[\lim_{u\to \infty}
\frac{[\mathfrak{L}'(u)]^{\frac{1}{2}}}
{u^{\frac{\rho-1}{2}}[{\boldsymbol L}_f(u)]^{\frac{1}{2}}\zeta(u)}=\frac{\rho}{2}
\] which, together with (\ref{rv2}), produces
\begin{equation} \label{et2}
\lim_{t\to 0^+}\frac{[\mathfrak{L}'(\Phi(t))]^\frac{1}{2} [\Phi(t)]^{\frac{-\rho+1}{2}}}
{[{\boldsymbol L}_f(\Phi(t))]^{\frac{1}{2}}\int_0^t \mathfrak{K}(s)\,ds}=
\frac{\rho}{2}. \end{equation}
By (\ref{et1}), (\ref{et2}) and L'Hospital's rule, we find
\begin{equation} \label{ch1}
\lim_{t\to 0^+}\frac{\log \Phi(t)}{\log (\int_0^t \mathfrak{K}(s)\,ds)}=
\lim_{t\to 0^+}\frac{\Phi'(t)}{\Phi(t)}\,\frac{\int_0^t \mathfrak{K}(s)\,ds}{\mathfrak{K}(t)}=-
\frac{2}{\rho}.
\end{equation}
We differentiate (\ref{et1}) to obtain
\begin{equation} \label{spo}
\begin{aligned}
\Phi''(t)=\frac{-\mathfrak{K}(t) \Phi'(t) [{\boldsymbol L}_f(\Phi(t))]^{\frac{1}{2}} }
{ [\mathfrak{L}'(\Phi(t))]^{\frac{1}{2}}[\Phi(t)]^{\frac{1-\rho}{2}}}
& \left\{  \frac{\rho+1}{2}+\frac{\mathfrak{K}'(t)\Phi(t)}{\mathfrak{K}(t)\Phi'(t)}+
\frac{\Phi(t){\boldsymbol L}_f'(\Phi(t))}{2{\boldsymbol L}_f(\Phi(t))} \right. \\
& \ \  \left.  -\frac{\Phi(t)\mathfrak{L}''(\Phi(t))}{2\mathfrak{L}'(\Phi(t))}
\right\}
\end{aligned}
\end{equation}
for each $t\in (0,\beta)$.
By $\mathfrak{K}\in NRV_\theta(0+)$ we mean
$\widetilde {\mathfrak{K}}(u)=\mathfrak{K}(1/u)\in NRV_{-\theta}$. Hence,
$\lim_{t\to 0^+}\frac{t\mathfrak{K}'(t)}{\mathfrak{K}(t)}=\theta$ and
$\lim_{t\to 0^+}\frac{\int_0^t \mathfrak{K}(s)\,ds}
{t\mathfrak{K}(t)}=\frac{1}{\theta+1}$.
This, combined with (\ref{ch1}), yields
\begin{equation} \label{ha}
\lim_{t\to 0^+}\frac{\mathfrak{K}'(t)}{\mathfrak{K}(t)}\,\frac{\Phi(t)}{\Phi'(t)}=-\frac{\rho 
\theta}
{2 (\theta+1)}\quad \mbox{and} \quad
\lim_{t\to 0^+}\frac{t\Phi'(t)}{\Phi(t)}=-\frac{2(\theta+1)}{\rho}.
\end{equation} Thus, $\Phi\in NRV_{-\frac{2(\theta+1)}{\rho}}(0+)$.
By (\ref{ha}) and L'Hospital's rule, we obtain
\begin{equation} \label{end2}
\lim_{t\to 0^+}\frac{\log \Phi(t)}{\log t} =\lim_{t\to 0^+}
\frac{t\Phi'(t)}{\Phi(t)}=-\frac{2}{\rho}\,(1+\theta).
\end{equation}
Proceeding by induction, we conclude (\ref{jar}). Since
${\boldsymbol L}_f$ is a normalised slowly varying function and
$\mathfrak{L}'\in NRV_{-1}$, we have
\begin{equation} \label{spo2}
\lim_{t\to 0^+}\frac{\Phi(t){\boldsymbol L}_f'
(\Phi(t))}{{\boldsymbol L}_f(\Phi(t))}=0\quad \mbox{and}\quad
\lim_{t\to 0^+}\frac{\Phi(t)\mathfrak{L}''(\Phi(t))}{\mathfrak{L}'(\Phi(t))}=-1. \end{equation}
By (\ref{spo}), (\ref{ha}) and (\ref{spo2}), we infer that
\[ \lim_{t\to 0^+}\frac{\Phi''(t)[\mathfrak{L}'(\Phi(t))]^{\frac{1}{2}}}
{\mathfrak{K}(t)\Phi'(t)[\Phi(t)]^{\frac{\rho-1}{2}}[{\boldsymbol L}_f(\Phi)]^{\frac{1}{2}}}=
-\left(1+\frac{\rho}{2(\theta+1)}\right).\]
Replacing $\mathfrak{K}(t)$ by its value in (\ref{et1}), we obtain the first assertion
of (\ref{ch2}). Moreover,
\begin{equation} \label{ram0}
\lim_{t\to 0^+}\frac{\log (-\Phi'(t))}{\log \Phi(t)}=\lim_{t\to 0^+}
\frac{\Phi''(t)\Phi(t)}{[\Phi'(t)]^2}=1+\frac{\rho}{2(1+\theta)}.
\end{equation}
Since
$\mathfrak{L}$ varies slowly at $\infty$ and
$\mathfrak{L}'\in RV_{-1}$, we use
Proposition~\ref{reg} (i) to obtain
\begin{equation} \label{ram}
\lim_{t\to 0^+}\frac{\log \mathfrak{L}(\Phi(t))}{\log \Phi(t)}=0\quad
\mbox{and}\quad \lim_{t\to 0^+}\frac{\log \mathfrak{L}'(\Phi(t))}{\log \Phi(t)}=-1.
\end{equation}
We notice that
\[
\log \left( \frac{ \mathfrak{L}(\Phi(t)) \Phi(t)}{\mathfrak{L}'(\Phi(t))
[\Phi'(t)]^2}\right)
= \log \Phi(t)\left[
1+\frac{\log \mathfrak{L}(\Phi(t))}{\log \Phi(t)}-
\frac{2\log |\Phi'(t)|}{\log \Phi(t)}-
\frac{\log \mathfrak{L}'(\Phi(t))}{\log \Phi(t)} \right]
\]
which, together with (\ref{ram0}) and (\ref{ram}), leads to
\[\lim_{t\to 0^+}
\log \left( \frac{ \mathfrak{L}(\Phi(t)) \Phi(t)}{\mathfrak{L}'(\Phi(t))
[\Phi'(t)]^2}\right)=-\infty.
\]
Thus the second claim of (\ref{ch2}) is proved.
\end{proof}

\section{Proof of Theorem~\ref{teo1}}

Let us first remark that the Keller--Osserman condition $(A_2)$ holds.
Indeed, by using Proposition~\ref{reg}, we arrive at
\[ \lim_{z\to \infty} \frac{f(z)}{z^p}=\lim_{u\to \infty}
\frac{f(\mathfrak{L}(u))}{u^\rho {\boldsymbol L}_f(u)}\,\frac{u^\rho {\boldsymbol L}_f(u)}
{[\mathfrak{L}(u)]^p}=\lim_{u\to \infty} \frac{u^\rho {\boldsymbol L}_f(u)}{[\mathfrak{L}(u)]^p}=
\infty,\ \ \forall p>1.\]

Thus,
Eq. ($P$) has at least a large solution when
$a<\lambda_{\infty,1}$ and no large solution provided that
$a\geq \lambda_{\infty,1}$ (see \cite[Theorem~1.1]{ccm}).

We now prove that, for $a<\lambda_{\infty,1}$ fixed,
every large solution of ($P$) exhibits the same asymptotic behaviour
near $\partial \Omega$, namely (\ref{ii}).
Set
$\vartheta^{\pm}=\left(\frac{\rho}{2(1+\theta)(1\mp 2\epsilon_0)}
\right)^{\frac{1}{\rho}}$, where $\epsilon_0\in (0,1/2)$ is arbitrary.
Let $\delta\in (0,\beta/2)$ be small such that
\begin{itemize}
\item[(i)] $d(x)$ is a $C^2$-function on the set $\{x\in {\mathbb R}^N:\ d(x)<2\delta\}$.
\item[(ii)] $\mathfrak{K}$ is nondecreasing on $(0,2\delta)$.
\item[(iii)] $(1-\epsilon_0)\mathfrak{K}^2(d(x))<b(x)<
(1+\epsilon_0)\mathfrak{K}^2(d(x))$, for all $x\in \Omega$
with $d(x)<2\delta$.
\item[(iv)] $\mathfrak{L}(\vartheta^\pm \Phi(2\delta))>0$.
\end{itemize}

Let $\sigma\in (0,\delta)$ be arbitrary. We define
$u^\pm_\sigma(x)=\mathfrak{L}(\vartheta^\pm \Phi(d(x)\mp \sigma))$, where
$d(x)\in (\sigma, 2\delta)$ (resp., $d(x)+\sigma<2\delta$) for
$u_\sigma^+(x)$ (resp., $u^-_\sigma(x)$). It follows that
\[ \begin{aligned}
 \Delta u_\sigma^\pm &= {\rm div}
 \left( \vartheta^\pm \mathfrak{L}'(\vartheta^\pm \Phi(d(x)\mp \sigma))
 \Phi'(d(x)\mp \sigma)\nabla d(x)\right)\\
 & = (\vartheta^\pm)^2 \mathfrak{L}''(\vartheta^\pm \Phi(d(x)\mp \sigma))
 [\Phi'(d(x)\mp \sigma)]^2|\nabla d(x)|^2\\
& \quad +\vartheta^\pm \mathfrak{L}'(\vartheta^\pm \Phi(d(x)\mp \sigma))
\Phi''(d(x)\mp \sigma)|\nabla d(x)|^2\\
& \quad +
\vartheta^\pm \mathfrak{L}'(\vartheta^\pm \Phi(d(x)\mp \sigma))
\Phi'(d(x)\mp \sigma)\Delta d(x).
\end{aligned} \]

In view of (i)--(iii), when $\sigma<d(x)<2\delta $ we obtain
(since $|\nabla d(x)|=1$)
\begin{equation}\label{sta}
\begin{aligned}
\Delta u_\sigma^++au_\sigma^+-b(x)f(u_\sigma^+) \leq
\frac{ \vartheta^+ \mathfrak{L}'(\vartheta^+ \Phi(d(x)-\sigma))[\Phi'(d(x)-\sigma)]^2}
{\Phi(d(x)-\sigma)}\times \\
\times \left\{
\frac{\Phi(d(x)-\sigma)}{\Phi'(d(x)-\sigma)}\,\Delta d
+{\mathcal E}^+(d(x)-\sigma)
\right\}
\end{aligned}
\end{equation}
respectively, when $d(x)+\sigma<2\delta$,
\begin{equation}\label{sta1}
\begin{aligned}
\Delta u_\sigma^-+au_\sigma^--b(x)f(u_\sigma^-) \geq \frac{
\vartheta^- \mathfrak{L}'(\vartheta^-
\Phi(d(x)+\sigma))[\Phi'(d(x)+\sigma)]^2}
{\Phi(d(x)+\sigma)}\times \\
\times \left\{
\frac{\Phi(d(x)+\sigma)}{\Phi'(d(x)+\sigma)}\,\Delta d
+{\mathcal E}^-(d(x)+\sigma)\right\}
\end{aligned}
\end{equation}
Here ${\mathcal E}^\pm $
are real functions defined on $(0,2\delta)$ as follows
\begin{equation} \label{ao}
 {\mathcal E}^\pm(t):=
\frac{\Phi''(t)\Phi(t)}{[\Phi'(t)]^2}+\frac{\vartheta^\pm
\Phi(t)\,\mathfrak{L}''(\vartheta^\pm \Phi(t))}
{\mathfrak{L}'(\vartheta^\pm \Phi(t))}+\frac{a\mathfrak{L}(\vartheta^\pm \Phi(t))}{\vartheta^\pm
\mathfrak{L}'(\vartheta^\pm \Phi(t)))}\frac{\Phi(t)}{[\Phi'(t)]^2}
-{\mathcal D}^\pm(t)
\end{equation}
where we denote
\[
{\mathcal D}^\pm (t) =(1\mp \epsilon_0)\,
\frac{\mathfrak{K}^2(t) f(\mathfrak{L}(\vartheta^\pm \Phi(t))}
{\vartheta^\pm \mathfrak{L}'(\vartheta^\pm \Phi(t))}\,\frac{\Phi(t)}{ [\Phi'(t)]^2}. \]
By virtue of (\ref{et1}), we may rewrite
${\mathcal D}^\pm (t)$ as
\[{\mathcal D}^\pm (t)
 =(1\mp \epsilon_0)
 (\vartheta^\pm)^\rho\,
 \frac{f(\mathfrak{L}(\vartheta^\pm \Phi(t)))}
 {(\vartheta^\pm \Phi(t))^\rho {\boldsymbol L}_f(\vartheta^\pm \Phi(t))}\,
\frac{{\boldsymbol L}_f(\vartheta^\pm \Phi(t))}
{{\boldsymbol L}_f(\Phi(t))}\,\frac{\mathfrak{L}'(\Phi(t))(\vartheta^\pm)^{-1}}
{\mathfrak{L}'(\vartheta^\pm \Phi(t))}.\]
It follows that
$\lim_{t\to 0^+} {\mathcal D}^\pm (t)=(1\mp \epsilon_0)
(\vartheta^\pm)^\rho$.
Note that $\lim_{u\to \infty}\frac{u\mathfrak{L}''(u)}{\mathfrak{L}'(u)}=-1$ (since
$\mathfrak{L}'\in NRV_{-1}$).
Moreover, by (\ref{ao}) and Lemma~\ref{pro}, we find
\[ \lim_{t\to 0^+}{\mathcal E}^\pm(t)=\frac{\rho}{2(1+\theta)}-
(1\mp \epsilon_0)\,(\vartheta^\pm)^\rho=
-\frac{\rho}{2(1+\theta)}\,\frac{\pm \epsilon_0}{1\mp 2\epsilon_0}.\]
Hence, using (\ref{sta}) and (\ref{sta1}),
we can choose $\delta>0$ small enough so that
\begin{equation}
\begin{split}
\Delta u_\sigma^++au_\sigma^+-b(x)f(u_\sigma^+) &\leq 0,\quad
\forall x\ \mbox{with}\ \sigma<d(x)<2\delta\\
\Delta u_\sigma^-+au_\sigma^--b(x)f(u_\sigma^-) &\geq 0,\quad
\forall x\ \mbox{with}\ d(x)+\sigma<2\delta.
\end{split}
\label{form}
\end{equation}
For $\eta>0$, we define $C_\eta=\{x\in \Omega:\ d(x)>\eta\}$ and
$ D_\eta=\{x\in {\mathbb R}^N:\ \mathrm{dist}(x,\Omega)<\eta\}$.

Let $\eta>0$ be small such that
$a<\lambda_1(-\Delta, D_\eta\setminus\overline{\Omega})$, where
$\lambda_1(-\Delta, D_\eta\setminus \overline{\Omega})$ denotes the first
Dirichlet eigenvalue of $(-\Delta)$ in the domain
$ D_\eta\setminus \overline{\Omega}$.
Let $p\in C^{0,\mu}(\overline{D}_{2\eta})$ satisfy $0<p(x)\leq b(x)$
for $x\in \Omega\setminus C_{2\delta} $, $p=0$ on $\overline{D}_{\eta}
\setminus \Omega$ and $p>0$ on $\overline{D}_{2\eta}\setminus
\overline{D}_\eta$.

We denote by $w$ a positive large solution of $\Delta w+aw=p(x)f(w)$ in
$D_{2\eta}\setminus \overline{C}_\delta $ (see \cite[Theorem~1.1]{ccm}).
Set $U:=u_a+w$ and $V_\sigma:=u_\sigma^++w$, where $u_a$ is an arbitrary
large solution of ($P$). It follows that
\[\Delta U+a U-b(x)f(U)\leq 0\ \mbox{in}\ \Omega\setminus
\overline{C}_\delta\ \mbox{and}\
\Delta V_\sigma+a V_\sigma-b(x)f(V_\sigma)\leq 0\ \mbox{in}\
C_\sigma\setminus \overline{C}_\delta. \]
Notice that $U|_{\partial \Omega}=\infty>u_\sigma^-|_{\partial\Omega}$,
$U|_{\partial C_\delta}=\infty>u_\sigma^-|_{\partial C_\delta}$, resp.,
$V_\sigma|_{\partial C_\sigma}=\infty>u_a|_{\partial C_\sigma}$ and
$V_\sigma|_{\partial C_\delta}=\infty>u_a|_{\partial C_\delta}$.
Thus, by \cite[Lemma 2.1]{ccm}, we deduce
$u_a+w\geq u_\sigma^-$ on $\Omega\setminus \overline{C}_\delta$
and $u_\sigma^++w\geq u_a$ on $C_\sigma\setminus \overline{C}_\delta$.
Letting $\sigma\to 0^+$, we find
\[ \mathfrak{L}(\vartheta^-\Phi(d))-w\leq u_a\leq
 \mathfrak{L}(\vartheta^+\Phi(d))+w\quad \mbox{on}\ \Omega\setminus
\overline{C}_{\delta}.\] Since $w$ is uniformly bounded on
$\partial \Omega$ and $\mathfrak{L}$ varies slowly at $\infty$, we
conclude (\ref{ii}).

Thus, $\lim_{d(x)\to 0}\frac{u_1(x)}{u_2(x)}=1$ for any two large
solutions $u_1$, $u_2$ of ($P$). From now on, we can use the same
line of reasoning as in the proof of \cite[Theorem~1]{cr} to
obtain $u_1\equiv u_2$ on $\Omega$. \qed

\section*{Acknowledgment}
F. C\^{\i}rstea wishes to thank Prof. N. Dancer for bringing this
problem into her attention. She also thanks Prof. N. S. Trudinger
for fruitful discussions, hospitality and support received during
a short visit to the ANU in January 2004.

\end{document}